\documentclass[12pt]{article}
\usepackage{graphicx,amsfonts,amsmath,hyperref,amscd,epsf}
\usepackage{amsmath,amsfonts,amssymb,latexsym}
\usepackage{multicol}
\usepackage{epsfig,url,xr,latexsym}

\newcommand{\vp}{\varphi}
\newcommand{\ti}{\tilde}
\newcommand{\wti}{\widetilde}

\newcommand{\p}{\partial}
\newcommand{\M}{\mathcal{M}}
\newcommand{\g}{\gamma}
\newcommand{\om}{\omega}
\newcommand{\si}{\sigma}

\newcommand{\Om}{\Omega}

\newcommand{\mref}[1]{(\ref{#1})}
\newcommand{\beq}{\begin{equation}}
\newcommand{\eeq}{\end{equation}}
\newcommand{\beqn}{\begin{eqnarray}}
\newcommand{\eeqn}{\end{eqnarray}}
\def\Ga{\Gamma}
\def\la{\lambda}

\def\si{\sigma}
\def\al{\alpha}
\def\a{\alpha}

\def\la{\lambda}
\def\ve{\varepsilon}
\def\be{\beta}

\def\ga{\gamma}

\def\ti{\widetilde}

\def\g#1#2{\frac{\partial\,#1}{\partial\,#2}}

\def\hat{\widehat}
\def\tilde{\widetilde}
\def \bfo {\begin {eqnarray*} }
\def \efo {\end {eqnarray*} }
\def \ba {\begin {eqnarray*} }
\def \ea {\end {eqnarray*} }
\def \beq {\begin {eqnarray}}
\def \eeq {\end {eqnarray}}

\def \det {\hbox{det}}

\def\s{{\sigma}}
\def \p {\partial}

\def\M{{\mathcal M}}

\def\g{{\gamma}}

\newtheorem{definition}{Definition}[section]
\newtheorem{theorem}[definition]{Theorem}
\newtheorem{lemma}[definition]{Lemma}

\newtheorem{proposition}[definition]{Proposition}
\newtheorem{corollary}[definition]{Corollary}
\newtheorem{remark}[definition]{Remark}
\newtheorem{condition}[definition]{Condition}

\title{Inverse Boundary Spectral Problem for Riemannian Polyhedra}
\author{Kirpichnikova A., Kurylev Ya.}

\begin{document}
\maketitle

{\small We consider an admissible Riemannian polyhedron with
piece-wise smooth boundary. The associated Laplace defines the
boundary spectral data as the set of eigenvalues and restrictions to
the boundary  of the corresponding eigenfunctions. In this paper we
prove that the boundary spectral data prescribed on an open subset
of the polyhedron boundary determine the admissible Riemannian
polyhedron uniquely.}

\section{Introduction}\label{intro}
Recent years have seen some very significant achievements in the
study of inverse boundary-value problems in a single component body.
Mathematically such body is described by a PDE or a system of PDE's
with relatively smooth coefficients. Starting from the pioneering
works \cite{Bel} and \cite{SU}, inverse boundary-value problems were
solved, at least on the level of uniqueness and sometimes
conditional stability, for a wide range of scalar inverse problems,
both isotropic and anisotropic, see e.g. \cite{AKKLT},
\cite{Belishev}, \cite{BelKur}, \cite{KK2}, \cite{KSU}, \cite{LTU},
\cite{LeU}, \cite{Na1}, \cite{Na2}, \cite{Nov}, \cite{PU},
\cite{Syl} for a far from complete list of references, with further
references in monographs \cite{Isakov1} or \cite{KKL}. Moreover, for
such media there appeared a number of important results in the study
of the inverse boundary-value problems for systems of PDE's
corresponding to physically important models of electromagnetism,
elasticity and Dirac equations, see e.g. \cite{Iso}, \cite{KLS},
\cite{KL-Dirac}, \cite{NakUhl}, \cite{Nak-Dirac}, \cite{OPS},
\cite{OS}.

Much less is known, however, about the inverse boundary-value
problems for a multicomponent medium. Mathematically, such medium is
described by PDE or system of PDE's with piece-wise smooth
coefficients with different subdomains of the regularity of
coefficients corresponding to different components of the medium.
Clearly, the study of inverse problems for the multicomponent media
is of substantial importance for practical applications. Imagine,
for example, a human body with bones, muscle tissue, lungs, etc.
each of those having distinctive values of material parameters, or
an upper crust of the Earth which is a composition of clay, sand,
rock, oil, water, etc. A complete answer to the inverse boundary
problems in a multicomponent medium, at least when the data are
measured on the whole boundary, is obtained only for the
two-dimensional case. Namely, it was shown in   \cite{PaiAst},
\cite{ALP} that the Calderon inverse boundary problem in the 2D case
has a unique solution in the class of $L^{\infty}-$coefficients.
Clearly, these results cover also the case of a multicomponent
medium. In higher dimensions, the results are restricted mainly to
the inverse obstacle problem.  In these problems the goal is to find
a shape of an inclusion inside a given medium which parameters are
known {\it a priori}.  In the case when parameters of a medium
and/or inclusion are unknown they are assumed to be homogeneous
throughout each component, see e.g. \cite{AMR}, \cite{Ike},
\cite{Ike1}, \cite{Isakov}, \cite{KirPai}. Having said so, we should
note that there exist powerful methods to find singularities for
coefficients of lower order, see e.g. \cite{GrUhl}.

This paper is devoted to the study of the inverse boundary spectral
problem for the Laplace operator in a multicomponent medium. To be
more precise, we assume that the domain occupied by the medium
consists of a finite number of subdomains with piece-wise smooth
boundaries between them. The metric tensor in each subdomain is
smooth  but does have jump singularity across the interfaces, i.e.
the boundaries between adjacent subdomains. Adding proper
transmission conditions across the interfaces and boundary
conditions on the domain's boundary, defines a Laplace operator
which, from the spectral point of view, has effectively the same
properties as the Laplace operator in a single component medium.
Mathematically, the considered medium may be described as a
Riemannian polyhedron. Leaving exact definitions of an appropriate
Riemannian polyhedron to the next section, imagine an
$n-$dimensional simplicial complex ${\mathcal M}$ where simpleces
can be glued together, pairwise, along their $(n-1)-$dimensional
faces which we continue to call interfaces (sometimes $(n-1)-$
interfaces). Imagine now that each simplex has its own smooth metric
$g$ which, in principle, may have  jumps across interfaces between
adjacent simpleces. This, together with some additional
geometric/combinatoric conditions described in section \ref{prelim},
defines a Riemannian polyhedron $({\mathcal M}, g)$. Starting from
the corresponding Dirichlet form on $H^1({\mathcal M}, g)-$functions
and using standard methods of spectral theory, the Laplace operator
with Neumann boundary conditions, $\Delta$,  is then well-defined in
$L^2({\mathcal M}, g)$. Denote by $\{\la_k,\, \varphi_k\}_{k=1}^{
\infty}$ the set of all eigenvalues, counting multiplicity, and
corresponding orthonormal eigenfunctions of $\Delta$. Let $\Gamma
\subset \p {\mathcal M}$ be open.

\begin{definition}
\label{bsd} The collection $\left(\Gamma,\,\{\la_k,\,
\varphi_k|_{\Gamma}\}_{k= 1}^{\infty}\right)$ is called the (local)
boundary spectral data (LBSD) of the Riemannian polyhedron
$({\mathcal M}, g)$.
\end{definition}

Let now $({\mathcal M}, g)$ and $({\tilde{\mathcal M}}, {\tilde g})$
be two Riemannian polyhedra with LBSD $\left(\Gamma,\, \{\la_k,\,
\varphi_k|_{\Gamma}\}_{k=1}^{\infty}\right)$ and $\left({\tilde
\Gamma},\, \{{\tilde \la}_k,\, {\tilde \vp}_k|_{{\tilde
\Gamma}}\}_{k=1}^{\infty}\right)$, correspondingly.

\begin{definition}
\label{equivalence} LBSD for $({\mathcal M}, g)$ and
$({\tilde{\mathcal M}}, {\tilde g})$ are equivalent if
\begin{enumerate}

\item
$\Gamma$ and ${\tilde \Gamma}$ are homeomorphic, $\varkappa: \Gamma
\rightarrow {\tilde \Gamma}$;

\item
$\la_k= {\tilde \la}_k, \quad k=1, 2, \dots;$

\item
If $\la_k$ has multiplicity $m+1,\,m=0,1,\dots$, i.e.
$\la_k=\la_{k+1}=\la_{k+m}$, then there is an $(m+1) \times (m+1)$
unitary matrix $\mathbb{U}_k$ such that \bfo (\vp_k|_{\Gamma},
\dots, \vp_{k+m}|_{\Gamma})= \mathbb{U}_k \,(\varkappa^* {\tilde
\vp}_k|_{{\tilde \Gamma}}, \dots, \varkappa^*{\tilde
\vp}_{k+m}|_{{\tilde \Gamma}} ). \efo
\end{enumerate}
\end{definition}

We can now formulate the main result of the paper:

\begin{theorem}
\label{main} Let $({\mathcal M}, g)$ and $({\tilde{\mathcal M}},
{\tilde g})$ be two admissible Riemannian polyhedra. Let, in
addition, the metric tensors $g$ and $\tilde g$ do have jumps across
all $(n-1)-$interfaces in ${\mathcal M}$ and ${\tilde{\mathcal M}}$,
correspondingly. Assume that LBSD $\left(\Gamma,\, \{\la_k,\,
\vp_k|_{\Gamma}\}_{k=1}^{\infty}\right)$ and $\left({\tilde
\Gamma},\, \{{\tilde \la}_k,\, {\tilde \vp}_k|_{{\tilde
\Gamma}}\}_{k=1}^{\infty}\right)$ are equivalent. Then $({\mathcal
M}, g)$ and $({\tilde{\mathcal M}}, {\tilde g})$ are isometric.
\end{theorem}

Let us make some comments on this theorem:
\begin{enumerate}

\item
If $\Omega_i$ is an $n-$dimensional simplex of ${\mathcal M}$ with a
smooth metric $g_i$, then $g_i$ determines an inner metric on any
$l$ dimensional, $l<n$, simplex of ${\mathcal M}$ which lies in
$\Omega_i$ (here and later we assume each simplex to be close). In
particular, any $(n-1)-$interface $\gamma$ of ${\mathcal M}$ belongs
to two adjacent $n-$simpleces, which we often denote in such case
$\Omega_-$ and $\Omega_+$, and therefore, has two different  metric
tensors $g_-|_{\gamma}$ and $g_+|_{\gamma}$. By a metric tensor $g$
having a jump singularity across $\gamma$ we mean that, for any $p
\in \gamma$, \bfo g_-(p) \neq g_+(p). \efo This assumption is of a
technical nature and, in section \ref{generalization} we will
significantly weaken it.

\item
As shown in section \ref{prelim}, any Riemannian polyhedron has a
natural structure of a metric space. The isometry of $({\mathcal M},
g)$ and $({\tilde{\mathcal M}}, {\tilde g})$ is understood with
respect to these metric structures.

\item
The boundary $\p {\mathcal M}$ of a Riemannian polyhedron $({\mathcal M}, g)$
is itself a $(n-1)-$dimensional Riemannian polyhedron, probably disconnected. As
$\Gamma,\, {\tilde \Gamma}$ are open subsets of $\p {\mathcal M}, \,\p {\tilde {\mathcal M}}$,
by reducing them if necessary we assume that $\Gamma$ and $ {\tilde \Gamma}$
are open subsets of some $(n-1)-$dimensional simplex of
$\p {\mathcal M}, \,\p {\tilde {\mathcal M}}$, correspondingly. In the future, we will
always assume this condition to be true.

\end{enumerate}

The plan of the paper is as follows: In section \ref{prelim} we
provide some preliminary material on geometry of Riemannian
polyhedra and properties of the Laplace operator on them. Section
\ref{GB} is devoted to the description and some properties of the
non-stationary Gaussian beams on a Riemannian polyhedron. We prove
Theorem \ref{main} in sections \ref{chamber1} and \ref{globalRP}.
The last section \ref{generalization} is devoted to some
generalizations and open questions.

\section{Preliminary constructions}\label{prelim}
\subsection{Admissible Riemannian polyhedron}
\begin{figure}
  \includegraphics[width=350pt]{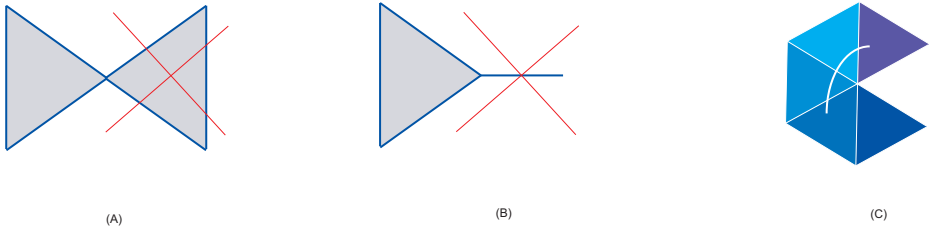}\\
  \caption{Case (A) is prohibited because its structure is not $(n-1)$-
  chainable; Case (B) is prohibited as it is not dimensionally
  homogeneous; Case (C) is appropriate}\label{Fig1}
\end{figure}

In this section we will introduce, following mainly \cite{FugEells}
and \cite{WBallman},  an admissible Riemannian polyhedron which is
the main object of the paper. We start with a closed $n-$dimensional
finite simplicial complex \bfo {\mathcal M} = \bigcup_{i=1}^I
\Omega_i, \efo where $\Omega_i$ are closed $n-$dimensional simpleces
of ${\mathcal M}$, with $\Omega_i^{\rm int}$ standing for the
interior of $\Omega_i$ which is an open subset of ${\mathcal M}$. We
assume that ${\mathcal M}$ is dimensionally homogeneous, i.e. any
$k-$simplex, $0\leq k \leq n$, of ${\mathcal M}$ is contained in at
least one $\Omega_i$.  We assume  also that any $(n-1)-$dimensional
simplex $\gamma$ belongs either to  two different $n$ simpleces,
$\Omega_i$ and $\Omega_j$, which in this case we will often denote
by $\Omega_-$ and $\Omega_+$, or to only one $n$ simplex $\Omega_i$.
In the former case we call $\gamma$ an {\it interface} (sometimes
$(n-1)-$dimensional interface) between $\Omega_-$ and $\Omega_+$, in
the latter case we call $\gamma$ a {\it boundary $(n-1)-$simplex}
with $(n-1)-$simpleces having this property forming the boundary $\p
{\mathcal M}$. We denote by ${\mathcal M}^k, \, 0\leq k \leq n$ the
$k-$skeleton of ${\mathcal M}$ which consists of all $k-$simpleces
of ${\mathcal M}$. Clearly, ${\mathcal M}={\mathcal M}^n$. We use
notations \bfo {\mathcal M}^{\rm int} = \bigcup \Omega_i^{\rm int},
\quad {\mathcal M}^{\rm reg} = {\mathcal M} \setminus \left(
\bigcup_{k=0}^{n-2} {\mathcal M}^{k}\right). \efo Following
\cite{FugEells}, we assume that ${\mathcal M}$ is $(n-1)-$chainable,
i.e. ${\mathcal M}^{\rm reg}$ is path connected, see Fig.
\ref{Fig1}.

Assume now that each $n-$simplex $\Omega_i$ is equipped with a
smooth (up to $\p \Omega_i$) Riemannian metric $g_i$, i.e.
$(\Omega_i,\,g_i)$ is a smooth Riemannian manifold with a piecewise
smooth boundary. This makes it possible to introduce the arclength
for admissible paths  $\eta:[0,a] \to {\mathcal M}$. We call a path
$\eta$ admissible if $\eta^{-1} ({\mathcal M}^{\rm int}) \subset
[0,a]$ is a (relatively) open subset of $[0,a]$ of full measure and,
if $\eta(\alpha, \beta)$ is in some $n-$simplex $ \Omega^{\rm int}$,
then $\eta: (\alpha, \beta) \to \Omega^{\rm int}$ is piecewise
smooth. Naturally, the arclength $|\eta(\alpha, \beta)|$ of the path
$\eta$ between $\eta(\alpha)$ and $\eta(\beta)$ is taken as \beq
\label{arclength} |\eta(\alpha, \beta)|= \int_\alpha^\beta \left[
g_{mj}(\eta(t))\dot\eta_m(t) \dot\eta_j(t)\right]^{1/2}dt, \eeq
where $\eta_j(t), \, \alpha <t <\beta$ are, for example, baricentric
coordinates in $\Omega$. As $\eta^{-1} ({\mathcal M}^{\rm int}) \cap
(0,a)$ consists of at most a countable number of open intervals
$(\alpha_i, \beta_i)$ we define \beq \label{length} |\,\eta[0,a]| =
\sum_i |\,\eta(\alpha_i, \beta_i)|. \eeq Next we introduce, for any
$p, q \in {\mathcal M}$, the distance, $d(p,q)$, \bfo d(p,q) =
\inf_{\eta}|\eta|, \efo where infenum is taken over all admissible
paths connecting $p$ and $q$. This makes $({\mathcal M},\,g)$ into a
metric space with its metric topology being the same as the topology
of a simplicial complex, see \cite{FugEells}.
\begin{definition}
\label{admissible} $({\mathcal M},\,g)$ is an admissible Riemannian
polyhedron if, for any $p, q \in {\mathcal M}$, \bfo d(p,q)
=\underset{\tilde \eta}{\inf}|\tilde \eta|, \efo where $\tilde \eta$
run over the subset of admissible paths between $p$ and $q$ such
that \bfo \tilde\eta^{\,-1} \left(\bigcup_{k=0}^{n-2} {\mathcal
M}^{k}\right) \setminus (\{0\}\cup\{1\})=\emptyset. \efo
\end{definition} As ${\mathcal M}$ is finite, the above condition is independent of a
particular choice of metric $g$.

\begin{figure}
  \includegraphics[width=380pt]{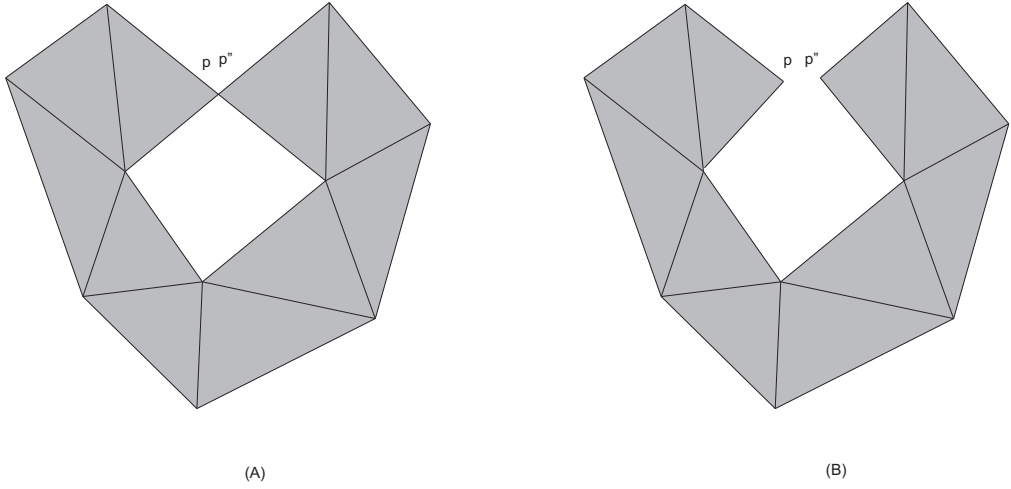}\\
  \caption{A sample of two topologically different Riemannian polyhedra having the same
  spectral properties.  It is stronger than the
$(n-1)-$chainability and is aimed at avoiding topologically
different Riemannian polyhedra which, however, can have the same
spectral properties}\label{Fig2}
\end{figure}

\subsection {Boundary normal and interface coordinates}\label{sec22}
In addition to the baricentric coordinates in any $\Omega^{\rm
int}$, we will often use boundary normal or interface coordinates
associated with $(n-1)-$subsimpleces of $\Omega$.

Let first $\gamma \in \Omega \cap \p {\mathcal M}$ be  a boundary
$(n-1)-$dimensional simplex with its $(n-1)-$dimensional interior
denoted by $\gamma^{\rm int}$. We introduce {\it boundary normal}
coordinates in an relatively open subset $\mathbf{U}\subset \Omega$
as \bfo p \to (s(p), \sigma(p)), \quad p \in U.\efo Here $\sigma(p)
= d(p, \p {\mathcal M})$, and we assume that there is a unique $q
\in \gamma^{\rm int}$ with $d(p, \p {\mathcal M})=d(p,q)$, such that
$p$ lies on the normal geodesic to $q$, $\varsigma_{\nu}(\tau), \,
\varsigma_{\nu}(0)=q,\, \varsigma_{\nu}(d(p, \p {\mathcal M}))=p$
and $\varsigma_{\nu}(0, d(p, \p {\mathcal M})) \in \Omega^{\rm
int}$. If $s(q)=(s^1,\dots,s^{(n-1)})$ are some (local) coordinates
on $\gamma$, e.g. baricentric coordinates, then $s(p)=s(q)$.

Let now $\g \subset \Omega_-\cap \Omega_+$ be an $(n-1)-$interface
between $n-$simpleces $\Omega_-$ and  $\Omega_+$. Let
$\mathbf{U}_{\pm} $ be relatively open subsets of $ \Omega_{\pm}$
with the nearest point on $\p \Omega_{\pm}$ lying on $\g$ such that
$(\mathbf{U}_- \bigcap \g) =(\mathbf{U}_+ \bigcap \g)$. Denote by
$(s, \s_{\pm})$ the boundary normal coordinates in $\mathbf{U}_{\pm}
$, where $s=(s^1,\dots, s^{(n-1)})$ are some local coordinates on
$\g$, e.g. baricentric coordinates with respect to $\Omega_-$ or
$\Omega_+$. We introduce the interface coordinates $(s, \s)$ on
$\mathbf{U}_- \bigcup \mathbf{U}_+$: \bfo (s,\s)=
\begin{cases}
(s,-\s),\,\,\hbox{in}\,\, \mathbf{U}_-,
\\ (s,\s),\,\, \hbox{in} \,\, \mathbf{U}_+.\,\,\end{cases}
\efo Then the metric element in these coordinates takes the form,
\beq\label{element} (dl)^2 =(d\s)^2+
(g_{\pm})_{\alpha\beta}(s,\s)\,ds^{\alpha} ds^{\beta}. \eeq
Throughout the paper we assume that the following condition takes
place:
\begin{condition}
\label{A}
For any interface $\g$ and any point
$q=(s,0)$ on $\g$, the metric tensor $g_{\a\beta}$ has a jump
singularity at $q$.
\end{condition}

\subsection{Laplace operator}\label{laplace}
Let $H^1({\mathcal M})$ be the Sobolev space of functions $u \in
L^2({\mathcal M})$ such that $u_i=u|_{\Omega_i} \in H^1(\Omega_i)$
and, for any interface $\g$ between $\Omega_i$ and $\Omega_j$, \bfo
u_i|_{\g}=u_j|_{\g}. \efo The inner product on $H^1({\mathcal M})$
determines the closed non-negative Dirichlet form, \bfo {\Bbb
D}[u,v]= \sum_{i=1}^I(u_i, v_i)_{H^1(\Omega_i)}, \quad u,v \in
H^1({\mathcal M}). \efo By the standard technique of the theory of
quadratic forms, the form $\Bbb D$ determines a self-adjoint
operator in $L^2({\mathcal M})$, namely, the Laplace operator with
Neumann boundary condition, $\Delta$. The domain ${\mathcal
D}(\Delta)$ is defined by \beq \label{domain} {\mathcal D}(\Delta)
=\{u \in H^1({\mathcal M}):\,\, {\Bbb D}[u,v]=(f,v)_{L^2({\mathcal
M})} \,\, \hbox{for some}\,\, f \in L^2({\mathcal M})\}, \eeq where
$v \in H^1({\mathcal M})$ is arbitrary. Analysing condition
(\ref{domain}), we see that $u \in {\mathcal D}(\Delta)$ if $u \in
H^2(\Omega_i),\, i=1, \dots, I$ and, on any interface $\g \subset
\Omega_- \bigcap \Omega_+$, \beq \label{interfacecondition}
u_-|_{\g}=u_+|_{\g}, \quad [\sqrt {g_-} \p_{\s}u_-]|_{\g}=[\sqrt
{g_+} \p_{\s}u_+]|_{\g}. \eeq where $g_{\pm}(s)=
\det[\,(g_{\pm})_{\a\beta}(s,0)]$.

As, due to the finiteness of ${\mathcal M}$, the embedding of
$H^1({\mathcal M})$ into $L^2({\mathcal M})$ is compact, the
spectrum of $\Delta$ is pure discrete, \bfo 0=\la_1 <\la_2 \leq
\dots, \quad \la_k \to \infty, \efo with the corresponding basis of
orthonormal eigenfunctions to be denoted by
$\{\vp_k\}_{k=1}^{\infty}$. Standard considerations, see e.g.
\cite{BelKur} or \cite{Taylor} show that $\{\vp_k\}_{k=1}^{\infty}$
distinguish points in ${\mathcal M}^{\rm int}$, i.e. for $p \neq q
\in {\mathcal M}^{\rm int}$, there is $k$ with $\vp_k(p) \neq
\vp_k(q)$. Moreover,
\begin{proposition}\label{disting} Let $\{\vp_k\}_{k=1}^\infty$ be
an orthonormal basis of eigenfunctions of the Laplace operator
$\Delta.$ Then $\{\vp_k\}_{k=1}^{\infty}$ form local coordinates
near any $p \in {\mathcal M}^{\rm int}$, i.e. there are
$k_1(p),\dots, k_n(p)$ such that $(\vp_{k_1},\dots, \vp_{k_n})$ form
local coordinates near $p$.\end{proposition}

\section{Gaussian Beams near interfaces}\label{GB}
\subsection{Gaussian beams on smooth manifolds}\label{GBsmooth}
In this section we briefly recall some results on the non-stationary
Gaussian beams on smooth manifolds. Their theory goes back to the
pioneering works \cite{BabichUlin}, \cite{Katchalov},
\cite{Ralston}. In our exposition we follow mainly section 2.4 of
\cite{KKL}. Non-stationary Gaussian beams are some (formal)
solutions of the wave equation \beq\label{waveEQ}U_{tt}-\Delta
U=0,\eeq which are concentrated, at each moment of time $t$, near a
point $x(t)$. The point $x(t)$ moves with a unit speed along a
geodesic on a smooth Riemannian manifold $(\mathcal N,\,h)$ with
$\Delta$ being the Laplacian corresponding to $(\mathcal N,\,h).$
Introducing a moving frame
$$y(t)=x-x(t),$$ a formal Gaussian beam has a form as a formal
series \beq\label{formalGB} U_\ve(t,y)\asymp M_\ve
\exp{\{-(i\ve)^{-1}\Theta(t,y)\}}\underset{l\ge 0}{\sum}u_l
(t,y)(i\ve)^l.\eeq Here $M_\ve=(\pi\ve)^{-\frac{n}4},\,\,0<\ve\ll1;$
$\Theta$ and $u_l,\,l=0,1,...,$ are formal series in powers of $y$.
They are usually represented as sums of homogeneous polynomials in
$y$ with coefficients depending on $t$,
$$
\Theta\asymp\underset{m\ge 1}{\sum}\theta_m(t,y),\,\quad
u_l=\underset{m\ge 1}{\sum}u_{lm}(t,y),
$$
$\theta_m$ and $u_{lm}$ being homogeneous polynomials on $y.$ The
polynomials $\theta_m$ and $u_{lm}$ are chosen so that, considered
as formal series with respect to $y$ and $(i\ve),$
\begin{equation}\label{wave}\partial_t^2 U_\ve-\Delta U_\ve=0.
\end{equation} Note that $"\asymp"$ exactly means that the formal
series \mref{formalGB} satisfies formally equation \mref{wave}.

The most important properties of the non-stationary Gaussian beams
are: \begin{description}
    \item[(a)] $\theta_1(t,y)=(\xi(t),y(t))=\xi_j(t)y^j(t)$, where
    $\xi_j(t)$ is the unit covector corresponding to the geodesic $x(t);$
    \item[(b)] $\theta_2(t,y)=\langle H(t)y,y\rangle$, where $H(t)$
    is a symmetric matrix, satisfying $Im \langle H(t)y,y\rangle\ge C(T) |y|^2$, for $-T<t<T.$
\end{description}\begin{remark}\label{constants}
From now on throughout this paper we use the following notations $C
$ (or, $C_L(t)$) is a generic constant, $C>0$, independent of $\ve$;
$\mu(L)$ is defined for sufficiently large positive integers $L$
such that $\mu(L)\rightarrow \infty$ when $L\rightarrow \infty$.
\end{remark}
Conditions (a) and (b) imply that $U$ decays exponentially outside
an $\ve^{1/2}$ - neighborhood of $x(t)$. It is important to note
that, starting from a formal Gaussian beam $U_\ve$ we can construct
a family of solutions to the wave equation (\ref{waveEQ}), which
"looks like" $U_\ve$. To this end, we start with a finite series
\beq\label{finite} U^L_{\ve}=M_\ve
\exp{\{-(i\ve)^{-1}\Theta^L(t,y)\}}\underset{l=0}{\overset{L}{\sum}}
u^L_l(t,y)\,\chi(d^2(x,x(t))\ve^{-5/6}), \eeq
$$
\Theta^L=\underset{l=0}{\overset{L}{\Sigma}}\theta_m;\quad
u^L_l=\underset{l=0}{\overset{L}{\Sigma}}u_{lm},
$$ where $\chi(s)$ is a smooth cut-off function equal to $1$ near $s=0.$
Then
$$
\|\p^2_t U_\ve^L-\Delta U_\ve^L\|_{C^{\mu(L)}(\mathcal
N\times[-T,T])}\le C_L(T)\,\ve^{-\mu(L)}.
$$
By standard hyperbolic estimates there exists a solution $\mathcal
U_\ve^L$ to (\ref{waveEQ}) such that \beq\label{real}
\|(U_\ve^L-\mathcal U_\ve^L)\|_{C^{N}(\mathcal N\times[-T,T])}\le
C_L(T) \ve^{-\mu(L)}.\eeq

Moreover, if we generate a wave inside $\mathcal N$ by a boundary
source \beq\label{boundary} U_\ve|_{\p\mathcal N\times[-T,T]}=
f_\ve(t,s). \eeq Let $f_\ve(t,s),\,\,s\in\p\mathcal N,\,\,t\in
[-T,T]$ be given by a formal expansion
$$
f_\ve(t,s)\asymp\exp\{-(i\ve)^{-1}\widehat{\Theta}(t,s)\}\underset
{l\ge0}{\sum}\, \widehat{u}_l(t,s)(i\ve)^l,
$$ where
$$
\widehat{\Theta}(t,s)\asymp -t+\xi_\al
s^\al+<\widehat{H}((s,t),(s,t))>+
\underset{m\ge2}{\sum}\,\widehat{\theta}_m (t, s);\quad
g^{\al\be}(0)\xi_\al\xi_\be<1,
$$
\beq\label{boundary1}
\widehat{u}_l(t,s)\asymp\underset{m\ge0}{\sum}\widehat{u}_{lm}
(t,s), \eeq with $\widehat{\theta}_m,\,\,\widehat{u}_{lm}$ being
homogeneous polynomials of degree $m$ with respect to $(t,s)$.
Assume that $Im\,\hat H>0.$ Then there is a unique formal Gaussian
beam $U_\ve$ satisfying (\ref{wave}), and the boundary condition
\mref{boundary} for $-t_0<t<t_0$ with some $t_0>0$ depending only on
geometry of $(\mathcal N, \p\mathcal N, h)$. Moreover, the
corresponding geodesic $x(t)$ starts, at $t=0,$ from the point
$s=0,$ into the (co)diversion $(\xi_\al,\xi_n)$ with
$\xi_n=[g^{\al\be}\xi_\al\xi_\be]^{1/2}$ geodesic starting at $s=0$.
We will refer to this result saying that we can guarantee a
non-stationary Gaussian beam propagating transversally to
$\p\mathcal N$ by a proper choice of a boundary source (for these
and other results on non-stationary Gaussian beams see e.g.
\cite{KK2}, \cite{KKL}).

\subsection{Gaussian beams at interfaces}\label{GBinterface}

In this section we consider reflection and transmission of the
non-stationary Gaussian beams from and through an
$(n-1)-$dimensional interface $\ga$ between two $n$-simpleces
$\Om_-$ and $\Om_+$. As our constructions will be of a local nature
we can, without loss of generality, restrict them to the
reflection/transmission of the Gaussian beams from a smooth
interface inside a smooth manifold. These questions, for the
incidence angle less then critical, were considered in detail in
\cite{KirpichGb} and \cite{PopovMM}, in the latter restricted to the
isotropic media. Assuming that $x(0)\in\ga$ and introducing the
interface normal coordinates $(s,\sigma)$ with $s=0$ corresponding
to $x(0),$ we have, for the incident Gaussian beam, $U_\ve^{in}$
\beq\label{incidentonboundary} U^{in}_\ve(t,s)|_{\sigma=0}\asymp
M_\ve \exp\{(i\ve)^{-1}\widehat{\Theta}^{in}(t,s)\}
\underset{l=0}{\overset{\infty}{\sum}}(i\ve)^l u_l^{in}(t,s) \eeq
$$
\left[\sqrt{g_-}\p_\sigma U_\ve^{in}(t,s)\right]|_{\sigma=0}\asymp
M_\ve \exp\{(i\ve)^{-1}\widehat{\Theta}^{in}(t,s)\}
\underset{l=-1}{\overset{\infty}{\sum}}(i\ve)\,\widehat{u}_l^{\,in}(t,s).
$$ Here $\widehat{\Theta}^{in},\,u_l^{in},\,\widehat{u}_l^{in}$
are sums of homogeneous polynomials with respect to $(t,s)$ with
$$
\widehat{\Theta}^{in} (t,s)=-t+\xi_\al^{in}s^\al+\hat
H((s,t),(s,t))+...\,\,;
$$
$$
u_{-1,0}^{in}=\xi_n^{in}u_{0,0}^{in},\,\,\,\xi_n^{in}>0,\,\,\,
g^{\al\be}\xi_\al^{in}\xi_\be^{in}+(\xi_n)^2=1,
$$ with $(\xi_\al^{in},\xi_n^{in})$ being the (co)direction of the
Gaussian beam at $t=0$. When
$g_+^{\al\be}(0)\xi_\al^{in}\xi_\be^{in}<1,$ it is possible to
construct two formal non-stationary Gaussian beams $U_\ve^r$ and
$U_\ve^{tr}$ in $\Om_-$ and $\Om_+,$ correspondingly such that
\beq\label{total}U_\ve=\begin{cases}U_\ve^{in}+U_\ve^r,
\quad&\text{in}\,\,\Om_-,\\U^{tr}_\ve\quad\quad
&\text{in}\,\,\Om_+\end{cases}\eeq satisfies the wave equation
(\ref{wave}) and transmission conditions (\ref{interfacecondition}).
To this end we use the technique briefly described in section
\ref{GBsmooth}, which reduces the problem to finding boundary
conditions of form \mref{boundary}, \mref{boundary1} for $U_\ve^r$
and $U_\ve^{tr}$ at $\sigma=0.$ In turn, this is possible utilizing
transmission condition \mref{interfacecondition} if
$g_+^{\al\be}(0)\xi_\al^{in}\xi_\be^{in}<1.$ Summarizing
considerations of \cite{KirpichGb}, we obtain the following result
\begin{lemma} \label{formal}Let $U_\ve^{in}$ be a formal
non-stationary Gaussian beam which hits the interface $\sigma=0$ at
$s=0,\,t=0$ with its (co)direction $(\xi_\al^{in},\xi_\be^{in})$
satisfying $g_+^{\al\be}(0)\xi_\al^{in}\xi_\be^{in}<1.$ Then there
are two formal Gaussian beams $U_\ve^r$ in $\Om_-$ and $U_\ve^{tr}$
in $\Om_+$ such that the total wave $U_\ve$ satisfies the
transmission condition \mref{interfacecondition}. The (co)directions
$(\xi_\al^r,\xi_n^r)$ and $(\xi_\al^{tr},\xi_n^{tr})$ of the
geodesics, corresponding to $U_\ve^r$ and $U_\ve^{tr}$ satisfy at
$t=0,\,s=0,\,\sigma=0$ the equation (Snell's Law):
$$\xi_\al^r=\xi_\al^{tr}=\xi_\al^{in},\,\,\xi_n^{r}=-\xi_n^{in},\,\,
\xi_n^{tr}=(1-g^{\al\be}_+(\cdot)\xi_\al^{tr}\xi^{tr}_\be)^{1/2}.$$
The main amplitude coefficients $u^r_{0,0}$ and $u_{0,0}^{tr}$ of
$U_\ve^r$ and $U_\ve^{tr}$ are related to $u_{0,0}^{in}$ at
$t=0,\,s=0,\,\sigma=0$ by
\beqn\label{amplitude1}u^{tr}_{0,0}=\frac{2\sqrt{g_-}\, \xi_n^{in}}
{\sqrt{g_-}\,\xi_n^{in}+\sqrt{g_+}\, \xi_n^{tr}}u^{in}_{0,0};
\,\,\,\,u_{0,0}^r= -\frac{\sqrt{g_+}\,\xi_n^{tr}-\sqrt{g_-}\,
\xi_n^{in}}{\sqrt{g_-}\, \xi_n^{in} + \sqrt{g_+}\, \xi_n^{tr}}
u_{0,0}^{in}, \eeqn where $g_{\pm}=\det^{-1} \left[
g_\pm^{\al\be}(0) \right].$
\end{lemma}

We note that transmission condition \mref{interfacecondition} for
$U_\ve^{in}+U_\ve^r,\,$ and $U_\ve^{tr}$ is understood in the formal
sense. Namely, $U_\ve^r$ and $U_\ve^{tr}$ may be expressed in the
form \mref{formalGB} with $U_\ve^{r, tr}|_{\si=0},\, \sqrt{g_\pm}
\,\p_\si U_\ve^{r, tr}|_{\si=0}$ having decomposition of form
\mref{incidentonboundary}. Then \mref{interfacecondition} means that
$$\widehat{\Theta}^{in}=\widehat{\Theta}^r=\widehat{\Theta}^{tr}$$
and
$$u_l^{in}+u_l^r=u_l^{tr},\,\,\,\widehat{u}_l^{in}+\widehat{u}_l^r=
\widehat{u}_l^{tr}$$ as polynomial with respect to $(t,s)$.

Observe that the condition
$\left[g_-^{\al\be}(0)\right]^{n-1}_{\al,\be=1}\ne\left[g_+^{\al\be}(0)\right]
^{n-1}_{\al,\be=1}$ implies that $u_{0,0}^r\ne 0$ for almost all
$(\xi_\al^{in},\xi_n^{in}).$

When dealing with an incoming non-formal Gaussian beam
\mref{finite}, \mref{real}, which we will denote by $U_\ve^{in,L}$,
similar to the above we find $U_\ve^{r,L}$ and $U_\ve^{tr,L}$ by
formulae \mref{finite} with $\Theta^{r,L},\,\,u_l^{r,L}$ and
$\Theta^{tr,L},\,\,u_l^{tr,L}$ instead of $\Theta^{in,L},
\,\,u_l^{in,L}$, correspondingly. Clearly, they give use to
approximate transmission conditions
\beqn\label{approx1}\|[(U_\ve^{in,L}+U_\ve^{r,L} -U_\ve^{tr,L})
]\|_{C^\mu(\ga\times(-t_0,t_0))}\le C(T)\,\ve^{-\mu(L)},
\\\label{approx2}
\|[(\sqrt{g_-}\p_\sigma(U_\ve^{in,L}+U_\ve^{r,L})) -\sqrt{g_+}
\p_\sigma U_\ve^{tr,L}]\|_{C^{\mu(L)} (\ga\times(-t_0,t_0))} \le
C(T)\,\ve^{-\mu(L)}.\eeqn Add to $U^{tr,L}_\ve$ a function
$\Psi_\ve^L(t,s,\sigma),$
$$\Psi_\ve^L(t,s,\sigma)=\chi(\sigma)\underset{k=0}{\overset{L}\sum}
\sigma^k \Psi_k(s,t), $$ where $\Psi_k(s,t)$ are chosen so that
$$\widetilde{U}_\ve^{\al,L}=U_\ve^{\al,L}+\Psi_\ve^L$$ satisfies
\beqn\label{newboundary}&\left[\Delta^p (U_\ve^{in,L}+
\widetilde{U}_\ve^{r, L}) \right] |_\ga& =\Delta^p\,
U_\ve^{tr,L}|_\ga ; \\
\label{newboundary2} &\left[\sqrt{g_-}\,\p_\sigma\left(
\Delta^p(U_{\ve}^{in,L}+U_\ve^{r,L})\right)
\right]|_\ga&=\left[\sqrt{g_+}\,\p_\sigma\left(\Delta^p\,
U_\ve^{tr,L} \right)\right]|_\ga,\eeqn for $0\le p\le [\frac L2]$.
Clearly, \mref{approx1}, \mref{approx2} imply that
$$\|\Psi_\ve^L(s,t,\sigma)\|_{C^{\mu(L)}(\Om_+\times(-t_0,t_0))}
\le C(T)\,\ve^{-\mu(L)}.$$ Together with \mref{newboundary},
\mref{newboundary2}, it follows from the wave equation \mref{wave}
that we can modify $U_\ve^{r,L}$ and $\widetilde{U}_\ve^{tr,L,}$
$$\widetilde U_\ve^{r,L}=U_\ve^{r,L}+\Phi_\ve^{r,L};\,\,\,
\widetilde U_\ve^{tr,L}= {U}_\ve^{tr,L}+\Phi_\ve^{tr,L},$$ with
$$\|\p_t^{\mu(L)}\Phi_\ve^L\|_{L^2(\M)}\le C_L(t_0)\,\ve^{-\mu(L)};$$
$$\|\Delta^{[\mu(L)/2]} \Phi_\ve^L\|_{L^2(\M)}\le C_L(t_0)\, \ve^{-\mu(L)};$$ so
for $-t_0\le t\le t_0$. Therefore,
$$u_\ve^L=\begin{cases}U_\ve^{in,L}+U_\ve^{r,L},\,\,&\sigma\le 0,\\
U^{tr,L}_\ve,&\sigma\ge 0\end{cases}$$ satisfies the wave equation
\mref{wave} and coincide with $U_\ve^{in,L}$ for negative $t$, more
precisely, for $t\le -C\ve^{5/12}.$ Here $\Phi_\ve^L$ is given by
$\Phi_\ve^{r,L} $ for $\sigma<0$ and $\Phi_\ve^{tr,L}$ for
$\sigma>0.$ \begin{remark}\label{reflected}When $\ga$ is
$(n-1)-$simplex in $\p\M,$ we can modify the previous construction
to find the Gaussian beam reflected from $\ga.$ The part of Lemma
\ref{formal} related to $U_\ve^r$ remains valid with formula for the
main term $u_{0,0}^r$ taking the form \beq\label{amplitude11}
u_{00}^r=-u_{00}^{in}.\eeq Thus, by the considerations similar to
the above it is possible to find solutions to the wave equation
\mref{wave} which satisfy the Dirichlet boundary condition and look
like a Gaussian beam.
\end{remark}

\section{First $n-$simplex}\label{chamber1}\subsection{}
In this section we start proving the uniqueness Theorem \ref{main}.
Recall that we are given diffeomorphic open subsets
$\Gamma\subset\gamma,\,\,\widetilde{\Ga}\subset\widetilde{\ga},$
where $\ga$ and $\widetilde{\ga}$ are boundary $(n-1)-$simpleces of
$n-$simpleces $\Om\subset\M,$
$\,\,\widetilde{\Om}\subset\widetilde{\M}.$ We assume, after proper
unitary transformations in finite-dimensional spaces corresponding
to eigenvalues of higher multiplicity, that
\beq\label{bsd2}\la_k=\widetilde{\la}_k,\quad\varphi_k|_\Ga=
\varkappa^*\,\widetilde{\varphi}_k|_{\widetilde{\Ga}},\quad
k=1,2,...\,\,,\eeq where $\varkappa:\Ga\rightarrow\widetilde{\Ga}$
is a diffeomorphism. Our goal in this section is to prove the
following result \begin{lemma}\label{first}Let
$\varkappa:\Ga_0\rightarrow\widetilde{\Ga}_0,\,\,\Ga_0\Subset\ga,\,\,
\widetilde\Ga_0\Subset\widetilde\ga$ be a diffeomorphism satisfying
conditions of Definition \ref{equivalence}. Then there is an
isometry $X:\Om\rightarrow\widetilde{\Om},$ such that
\beq\label{local}\varphi_k|_\Om=X^*\,\widetilde{
\varphi_k}|_{\widetilde{\Om}};\quad k=1,2,...,\quad
X|_{\Ga_0}=\varkappa.\eeq\end{lemma} To prove this lemma, observe
first that if
$\Ga_0\Subset\Ga,\,\,\widetilde{{\Ga}}_0\Subset\widetilde{\Ga}$ and
$\tau>0$ satisfy:
$$\exp_{\p\Om}:\,\Ga_0\times\,\left[0,\tau\right)\rightarrow\M,\,\,
\exp_{\p\widetilde{\Om}}:\,\widetilde{\Ga}_0\times\,
\left[0,\tau\right)\rightarrow\widetilde{\M}$$ are regular, i.e. map
into $\Om$ and $\widetilde{\Om}$ correspondingly and provide an
diffeomorphism between $\Ga_0\times\left[0,\tau\right)$ and
$\exp_{\p\Om}(\Ga_0\times[0,\tau))\subset\Om$ and
$\widetilde{\Ga}_0\times[0,\tau)$ and $\exp_{\p\Om}(\widetilde{\Ga}
_0 \times[0,\tau)) \subset\widetilde{\Om}$ and, in addition,
$$d(\Ga_0,\M\setminus\Om),\,
d(\widetilde{\Ga}_0,\widetilde{\M}\setminus\widetilde{\Om})
>\tau,$$ then the diffeomorphism
$$\widehat{X}=\varkappa\times\mathbb{I}:\,\, \Ga_0\times(0,\tau) \rightarrow
\widetilde{\Ga}_0\times(0,\tau)$$ satisfies
$$\widehat{X}^*\,\widetilde{\varphi_k}|_ {\widetilde{\Ga}_0
\times(0,\tau)} = \varphi_k|_ {\Ga_0\times(0,\tau)}.$$ Moreover,
$\widehat X$ is actually an isometry between $\Ga_0\times(0,\tau)$
and $\widetilde{\Ga}_0\times(0,\tau)$ considered as domains in
$(\Om^{int},\,g)$ and $(\widetilde{\Om}^{int},\,\tilde g),$
correspondingly. The proof of this fact is identical to the smooth
case and is given in Section 4.4 of \cite{KKL}.

Assume now that $\om_1,\,\om_2\subset\Om^{int}$ and
$\widetilde{\om}_1,\,\widetilde{\om}_2\subset\widetilde{\Om}^{int}$
are open subsets with
$\widehat{X}_i:\,\om_i\rightarrow\widetilde{\om_i}$ being isometries
satisfying \beq\label{identifice}\widehat{X}_i^*\,
\widetilde{\varphi_k}|_{\widetilde{\om}_i}=\varphi_k|_{\om_i},\quad
i=1,2,\,\,k=1,2,...\,\,.\eeq By Proposition 1, $\hat X_1|_{\om_1
\cap \om_2}=\hat X_2|_{\om_1\cap\om_2} $ which makes it possible to
extend $\hat X_i$ into an isometry $\hat X:\om_1 \cup\om_2
\rightarrow\tilde \om_1\bigcup\tilde\omega_2$ which also satisfies
\mref{identifice}.

Consider the family of all pairs of open subsets $\om\subset\Omega$
and $\tilde\om\subset\widetilde\Om$ which are isometric to each
other and satisfy \mref{identifice}. Clearly, this family is
partially ordered by induction, by the above we can consider its
maximal element which we denote by $(\Om_m,\widetilde\Om_m)$ with
the corresponding isometry denoted by $\hat X_m.$ We want to show
that $\Om_m=\Om^{int},\,\,\widetilde\Om_m=\widetilde\Om^{int}.$

To proceed, recall the following result from \cite{KKL}, which is
proven for smooth manifolds but remains valid for Riemannian
polyhedra under the conditions formulated below.
\begin{theorem}\label{subdomains}\begin{enumerate}
\item Let $S\subset\Om$ (or $\tilde S\subset\widetilde\Om$) be a smooth
subdomain such that $\{\varphi_k(p)\}_k^\infty$ (or
$\{\tilde\varphi_k(\tilde p)\}_k^\infty$) are known for $p\in S$ (or
for $\tilde p\in \tilde S$). Assume that for $\tau>0$
$$exp_{\p S}:\p S\times (0,\tau)\rightarrow\Om\setminus S,\,\,exp_{\p\tilde
S}: \p\tilde S\times(0,\tau)\rightarrow\widetilde\Om\setminus\tilde
S$$ are regular. If, in addition, $\tau<\{d(S,\M\setminus\Om),\,
d(\tilde S,\widetilde\M\setminus\widetilde\Om )\}$ then these data
determine uniquely $\varphi_k|_{S_\tau},\,\tilde\varphi_k|_{\tilde
S_{\tau}},$ where $$S_\tau=S\bigcup {exp}_{\p S} (\p
S\times(0,\tau)),\,\,\tilde S_\tau=\tilde S \bigcup exp_{\p\tilde
S}(\p\tilde S\times(0,\tau)).$$
\item Let $u,\,\tilde u$ be solutions of the initial
boundary value problem, \beq\label{IBVPS}u_{tt}-\Delta u=F\in
C_0^\infty(S\times\mathbb{R}_+);\,\,\,\, u_{tt}-\tilde\Delta
u=\widetilde F\in C_0^\infty(\tilde S\times\mathbb{R_+}),\\\nonumber
u|_{t=0}=f\in C^\infty_0(S);\,\, \tilde u|_{t=0}=\ti f\in
C_0^\infty(\tilde S),\\\nonumber u_t|_{t=0}=\phi\in
C^\infty_0(S);\,\, \tilde u_t|_{t=0}=\ti \phi\in C_0^\infty(\tilde
S).\eeq Then these data determine $u,\,\,\tilde u$ on
$S_\tau\times\mathbb{R}_+,\,\tilde S\times\mathbb{R}_+,$
correspondingly. In particular, if $S$ and $\tilde S$ are isometric,
with isometry $X$ satisfying \mref{identifice}, there is an extended
isometry $X_\tau,$
$$X_\tau: S_\tau\rightarrow\tilde S_\tau$$ with \beq\label{new}
\varphi_k|_{S_\tau}= X_\tau^*\, \tilde\varphi_k|_{S_\tau},
\,\,u^f|_{S_\tau\times\mathbb{R}_+}= X_\tau^*\,\tilde u^{\tilde
f}|_{\tilde S_\tau\times\mathbb{R}_+}, \eeq when $f=X^*\,\tilde f,
\,\phi=X^*\,\tilde \phi,\, F=X^*\,\tilde F$.
\end{enumerate}\end{theorem}

\subsection{}\label{41} Based on Theorem \ref{subdomains}, we
will finish the proof of Lemma \ref{first}. 
Assume, that a maximal element $(\Om_m,\,\, \widetilde\Om_m)
\ne(\Om^{int},\, \widetilde\Om^{int}),$ where without loss of
generality we can take $\Om_m\ne\Om^{int}.$ Therefore, there is a
point $p\in \mathcal{C}l(\Om_m)\,\cap\,\Om^{int}.$ Observe, that as
$\hat X_m: \Om_m \rightarrow\widetilde \Om_m$ is an isometry, it may
be extended to the mapping $\hat X_m:\mathcal Cl(\Om_m) \rightarrow
\mathcal Cl ( \widetilde\Om_m).$ Consider the following possible
scenarios:
\begin{enumerate}
    \item  $\tilde p=\hat X_m(p)\in\widetilde\Om^{int}$. Denote by
    $\delta=\min (d(p,\p\Om), d(\tilde p,\p\widetilde\Om))$ and by
    $\rho =\min(i(\Om,g),\,i(\widetilde\Om,\tilde g)),$ where $i(N,h)$
    stands for the injectivity radius of the normal coordinates of
    Riemannian manifold $(N,h).$ Let $\delta_0=\frac14 \min(\delta,
    \rho)$ and $p_0\in\Om_m,\,\tilde p_0=\hat X_m(p_0)\in\widetilde
    \Om_m$ satisfy $d(\tilde p_0,\tilde p)<\delta_0.$ Let $0\le\si
    <\delta_0$ satisfies $B(p_0,\si)\subset\Om_m,\, B(\tilde p_0,\si)=
    \hat X_m(B(p_0,\si))\subset\widetilde\Om_m,$ where $B(p,r)$ is a
    closed ball of radius $r$ centered at $p.$ Taking $S,\,\tilde S$ in
    Theorem \ref{subdomains} to be $B(p_0,\si),\, B(\tilde p_0,\si)$,
    we see that conditions of this Theorem are satisfied for $\tau=2
    \delta_0$ with $S_\si=B(p_0,\si+2\delta_0),\,\widetilde
    S_\si=B(\tilde p_0,\si+2\delta_0)$. Therefore, $\hat X_m$ can be
    extended to $\Om_m\cup B(p_0, \si+2\delta_0)$ containing $p$,
    which contradicts the definition of $\Om_m.$
    \item  $\tilde p=\hat X_m(p)\in\ti\ga^{int},$ where $\tilde\ga$ is
    some $(n-1)-$subsimplex of $\widetilde\Om.$ Let now $\delta= \min
    (d(p,\p\Om)),\,d(\tilde p,\p(\widetilde\Om\cup\widetilde\Om_1)
    )$, where $\widetilde\Om_1$ is another $n-$simplex adjacent to
    $\tilde\ga$ (if $\tilde\ga\in\p\widetilde\M$ we take $\p(\widetilde
    \M\setminus \widetilde\Om)$ rather than $\p(\widetilde\Om \cup
    \widetilde\Om_1)$). As earlier, let $\delta_0=\frac14 \min (\delta,
    \rho, i_{\ti\ga}(\ti p)),$ where $i_{\ti\ga}(\ti p)$ is the radius
    of injectivity of the interface normal coordinates related to $\ti
    p,$ (or to the boundary normal coordinates of $\ti p\in\p\M$).
    Introduce $p_0,\,\ti p_0=\hat X_m(p_0)$ as in the case 1 and
    take balls $B(p_0,\si),\,B(\ti p_0,\si),\,0<\si<\delta_0$
    similar to the case 1. Consider now the non-stationary Gaussian
    beams $\widetilde U_\ve^{in,L}$ on $\widetilde\M$ which start
    at $t=0$ at $\tilde p_0,$ in direction close to the normal direction from $\tilde
    p_0$ to $\ti\ga,$ i.e. with the initial co-vector
    $\tilde\xi$ close to $(0,0,...,1).$ These Gaussian beams reflect
    from $\tilde\ga$ and return to $B(\tilde p_0,\si)$ approximately
    at the time $t=2\, d(\tilde p_0,\tilde\ga)-\si.$ Thus, for
    $\tilde\xi$ close to $(0,...,0,1),$ the total Gaussian beam
    $\widetilde U_\ve^L=\widetilde U_\ve^{r,L}+\wti U_\ve^{in,L}$ in
    $\widetilde\Om$ satisfies \beq\label{furthoroL}
    \max{|\widetilde U_\ve^L (\tilde q,t)|>C},\eeq when $\tilde q
    \in B(\tilde p_0,\si),\, t\in [2 d(\tilde p_0,\tilde\ga)-\si,\,
    2\, d(\tilde p_0,\tilde\ga)+\si]$.

    On the other hand, the corresponding Gaussian beam
    $U_\varepsilon^{in,L}$ in $\Om$ moves from $B(p_0,\si)$ and has
    no reflected part for $t\in(0,\,2\, d(\tilde p_0,\ti\ga)+
    \si))$, so that \beq\label{24a}\max{|U_\ve^L(q,t)|} \le C\,\ve^{-
    \mu(L)}, \quad\text{for}\,\,q\in B(p_0,\si). \eeq When $\ve$ is
    sufficiently small and $L$ is sufficiently large, \mref{furthoroL},
    \mref{24a} contradict the second equation \mref{new}, of Theorem
    \ref{subdomains}. As $\varphi_k|_{B(p_0,\si)}$ and $
    \ti\varphi_k|_{B(\ti p_0,\si)}$ are known we can evaluate the
    Fourier coefficients $u_k(t)=\ti u_k(t)$ of $U_\ve^L(t)$ and
    $\widetilde U_\ve^L(t)$ and, next, $$U_\ve^L(q,t)=\sum u_k(t)
    \varphi_k(q);\,\,\wti U_\ve^L(\ti q,t)=\sum u_k(t)
    \ti\varphi_k(\ti q),$$ for $\,q\in B(p_0,\si),\,\,\ti
    q_0\in B(\ti p_0,\si)$.
    \item  Let now $\tilde p\in\widetilde\M^{n-2}.$ Using e.g.
    baricentric coordinates in $\widetilde\Om,$ we see that there
    are $C>1$ and $\ti\delta>0$ such that if $\ti d
    (\ti p_0,\p\ti\Om)<\ti\delta,$ then there is a curve $\ti x:
    [0,1]\rightarrow\widetilde\Om^{int},\,$ $\ti x(0)=\ti p_0,\,
    \ti x(1)=\ti q_0$ such that:\begin{description}
        \item[(i)] $ d(\ti x(t),\p\widetilde\Om)> C^{-1} d(\ti
        p_0,\p\widetilde\Om)$;
        \item[(ii)] $|\ti x[0,1]|\le C\, d(p_0,\p\widetilde\Om)$,
        where $|\ti x[a,b]|$ is the arclength of $\ti x(t)$ between
        $\ti x(a)$ and $\ti x(b)$;
        \item[(iii)] $d(\ti q_0,\p\widetilde\Om)= d
        (\ti p_0,\p\widetilde\Om) $. However, we can assume that
        there is a unique nearest point $\ti q\in\p\widetilde\Om$ to
        $\ti q_0$ and, in addition,
        \item[(iv)] $\ti q_0\in\ti\ga^{int},$ where $\ti\ga$ is an
        $(n-1)-$interface between $\widetilde\Om$ and some
        $\widetilde\Om_1$ (or $\ti\ga\subset\p\widetilde\M$);
        \item[(v)] Interface (boundary) normal coordinates centered
        at $\ti q_0$ are regular in $4\, d(\ti q_0,\ti q)$-vicinity
        of $\ti q_0$.
    \end{description}

Returning to the consideration of the case $\ti p \in \widetilde
\M^{n-2}, $ let $p_0\in\Om_m$ satisfy $d(p_0,p)\le \min (\ti\delta,
\frac1{8C}\delta),$ where $\delta$ is defined as in the case 1. With
$\ti p_0=\hat X_m(p_0),$ let $\ti x(t)$ be a curve in $\widetilde
\Om^{int}$ described earlier. Then the diffeomorphism $\hat
X^{-1}_m$ can be extended onto some open neighborhood of $\ti x(t)$.
Indeed, by the construction of the step 1, we can move recurrently
along $\ti x(t),$ using the balls of radius $\frac1{2C} d (\ti
p_0,\p\widetilde\Om).$ By the maximality of $(\Om_m, \widetilde
\Om_m)$, $\ti q_0\in\widetilde\Om_m$ with $q_0=\hat X_m^{-1}(\ti
q_0)\in\Om_m$ satisfying \beq\label{3ineq}d(q_0,\p\Om)\ge 4\, d (\ti
q_0,\ti q). \eeq Inequality \mref{3ineq} makes it possible to use
the same considerations as in the case 2, proving that the case $\ti
p\in\widetilde\M^{n-2}$ is possible.
\end{enumerate} We finish the section by the following
\begin{corollary}\label{boundary?}
Let $X:\Om\rightarrow\wti\Om,$ where $\Om,\,\,\wti\Om$ are
$n$-simpleces of $\M$ and $\wti\M,$ satisfy \mref{local}. If
$\ga,\,\ti\ga=X(\ga)$ are $(n-1)-$simpleces of $\Om$ and $\wti\Om$
correspondingly. Then $\ga$ is an interface if and only if $\ti\ga$
is an interface.\end{corollary}\textbf{Proof} Let $\ti f=\ti
f_\ve^L, \,\ti \phi=\ti \phi_\ve^L$ be the initial data for Gaussian
beam $\ti U_\ve^{in,L}$ starting in $\widetilde\Om$ towards $\ti\ga$
and $(f_\ve^L,\, \phi_\ve^L)=X^*\, (\ti f_\ve^L,\,\ti \phi^L_\ve)$
be the initial data of the Gaussian beam $ U_\ve^{in,L}.$ Comparing
$\wti U^{r,L}_\ve$ and $U_\ve^{r,L}=X^*\, \wti U^{r,L}_\ve$ and
using formulae \mref{amplitude1} and \mref{amplitude11}, we obtain
the desired result. $\Box$
\section{Global Isometry}\label{globalRP}
\subsection{}\label{51}In this section we show that
isometry on $X:\Om\rightarrow\wti\Om$ which satisfies \mref{local}
can be extended to a global isometry $X:\M\rightarrow\wti\M$
satisfying \beqn\label{global}\varphi_k(p)=X^*\,\ti\varphi_k(\ti
p),\,\, k=1,2,...,\,\,p\in\M^{int},\,\ti p=X(p)\in\wti\M^{int},
\\\nonumber d(p,q)=d(\ti p,\ti q),\,\,\ti p=X(p),\,\ti q=X(q),
\,\,p,q\in\M.\eeqn We start with the following result which is a
partial generalization of Theorem \ref{subdomains}.
\begin{lemma}\label{crossing}Let $\ga$ be the interface between
$\Om_-$ and $\Om_+$ in $\M$ and $\ti\ga=X_-(\ga)$ is the interface
between $\wti\Om_-$ and $\wti\Om_+$ in $\wti\M.$ Assume that
$X_-:\Om_- \rightarrow\wti\Om_-$ is an isometry satisfying
\mref{local}. Then $X$ can be extended to an isometry $X:
\Om_-\cup\Om_+ \rightarrow\wti\Om_-\cup\wti\Om_+$ which satisfies
\mref{local} on $\Om_-^{int}\cup\Om_+^{int}.$
\end{lemma}\textbf{Proof} Let $D,\,\wti D$ be smooth subdomains in
$\Om_-^{int}$ and $\wti\Om_-^{int}$ with the "upper" part of the
boundary parallel to $\ga,\,\ti\ga,$ i.e. the parts of $\p
D,\,\p\wti D$ are given by $x=(s,-\si_0),\,\ti p=(\ti s,-\ti\si_0),
\,\, s\in\Ga \Subset \ga,\,\ti s\in\ti\Ga_0\Subset\ti\ga.$ Assume,
without loss of generality that $d(D,\p(\Om_-\cup\Om_+)),\,\, d(\wti
D,\p(\wti\Om_- \cup\wti \Om_+))> 10\si_0$, $\,\si_0< \frac18 (
\min({i_{\Ga},i_{\wti\Ga}})),$ where $i_\Ga,\,i_{\wti\Ga}$ are the
injectivity radii of interface coordinates related to
$\Ga,\,\wti\Ga,$ correspondingly. We want to show that equation
\mref{local} implies that \beq\label{contin}
\varphi_k(s,\si)=\ti\varphi_k(\ti s,\ti\si),\,\,\text{for }
\,\,s\in\Ga,\,\ti s=X_-|_{\ga}(s),\,\,-2\si_0 <\si <2\si_0.\eeq We
have the following rather straightforward generalization of Tataru's
approximate controllability: \beq\label{controllability} \mathcal
Cl_{L^2(\M)}\{ u^f(t_0),f\in C^\infty_0(D\times(0,t_0))\}= L^2
(U_{t_0}(D)),\,\,0\le t_0\le 4\si_0,\eeq with a similar identity for
$\wti D.$ Here $U_{t_0}(D)$ is a $t_0-$neighborhood of $D$ in $\M$
and $u^f$ is a solution to the initial boundary value problem
$$u^f_{tt}-\Delta u^f =f,\quad u^f|_{t<0}=0,\quad
u^f|_{\p\M\times\mathbb{R}_+}=0.$$ As usual this result follows
immediately from observability, see e.g. Section 2.5 of \cite{KKL}.
To formulate the desired observability, let \beqn\nonumber
\begin{cases} v_{tt}-\Delta v=0,\quad \text{in} \quad\M
\times(-t_0,t_0),\,\, v|_{\p\M\times(-t_0,t_0)}=0, \\
v|_{D\times(-t_0,t_0)}=0, \end{cases}\eeqn where $0<t_0<4\si_0$.
Then $v=0$ in double cone
$$K(D,t_0)=\{(p,t) \in \M \times(-t_0,t_0):\,\, d(p,D)<t_0-|t|\}.$$
To use this we first observe that for any $q$ with $d(p,D)<4\si_0$
and any $\delta>0$ there is a piece-wise smooth curve $x: [0,1]
\rightarrow \Om_-\cup\Om_+$ with $x(0)=p\in D,\,\,x(1)=q$ and
$|x[0,1]|<d(q,D)+\delta.$ Moreover, $x(t)$ may be chosen to cross
$\ga$ transversally. Then we can continue $v$ by $0$ along this
curve $x(t)$ so that for $x(s)$ $$v=0 \quad \text{in}\,\,V_\si\times
(-(t_0-|x[0,s]|,\,t_0- |x(0,s)|)),$$ where $V_s$ is a small vicinity
of $x(s)$. Indeed, this is obvious for pieces of the path $x(s)$
lying inside either $\Om^{int}_-$ or $\Om^{int}_+.$ To cross $\ga$
we just observe that, if $v,\,\p_\si v=0$ on $\Sigma\times(-\hat
t,\hat t\,),\,\Sigma\Subset\ga,\,$ when approaching $\ga$ from
$\Om_+$, then by (8), $$v,\,\p_\si
v=0,\quad\text{on}\,\,\Sigma\times(-\hat t,\hat t\,),$$ when
approaching $\ga$ from $\Om_-,$ and, therefore, $v$ can be continued
by $0$ further along $x(t)$ into $\Om_+$ and vice versa from $\Om_+$
into $\Om_-$. Clearly, to use \mref{interfacecondition} we should
assume $v$ to be sufficiently regular, however, by smoothing $v$
with respect to time $t,$ we can extend it to non-smooth solutions,
see e.g. \cite{KKL2}.

This implies that $v(q,t)=0$ for $|t|<t_0-d(q,D)-\delta$. As
$\delta>0$ is arbitrary, we obtain that $v=0$ in $K(D,t_0)$.

Identity \mref{controllability} makes it possible, starting from
$\{\lambda_k,\,\varphi_k|_{D}\}_{k=1}^\infty$, and from
$\{\ti\lambda_k,\,\ti\varphi_k|_{\wti D}\}_{k=1}^\infty$, such that
$\lambda_k=\ti\lambda_k$ and $\varphi_k(s,\si)=\ti\varphi_k(\ti s,
\si)$, for $-2\si_0<\si<-\si_0,\,\,\ti s=(X_-(\ga)(s)),\, s\in\Ga$
to construct $\varphi_k(s,\si)$ and $\ti\varphi_k(\ti s,\si)$ for
$2\si_0<\si<2\si,\,\,\ti s=(X_-(\ga))(s),\,s\in\Ga$, see \cite{KK2},
Chapter 4.4 of \cite{KKL}. In particular, the construction in
\cite{KK2}, \cite{KKL} imply that $\ti\varphi_k(X_-(\ga)(s),\si) =
\varphi_k (s,\si)$. Observe now that $\Ga\times[0,2\si_0],\, \ti\Ga
\times [0,2\si_0]$ from a relatively open subdomains in
$\Om_+,\,\wti\Om_+,$ respectively with $X_+:\,(s,\si) \rightarrow
((X_-|_{\ga})(s),\si),\,0\le\si\le 2\si,$ being a diffeomorphism
satisfying \mref{identifice}.

Mimicking the proof of Lemma \ref{first} in section \ref{41}, we
extend $X_+$ to be a diffeomorphism, $X_+:\,\Omega_+\rightarrow \wti
\Om_+,$ satisfying \mref{local}. As, by construction, $X_+|_\ga =
X_-|_\ga$, $X$ defined as $X_-$ on $\Om_-$ and $X_+$ on $\Om_+$ is a
desired isometry between $\Om_-\cup\Om_+$ and
$\wti\Om_-\cup\wti\Om_+$.

\subsection{Identification of $n$-simpleces} 
\begin{figure}
  \includegraphics[width=360pt]{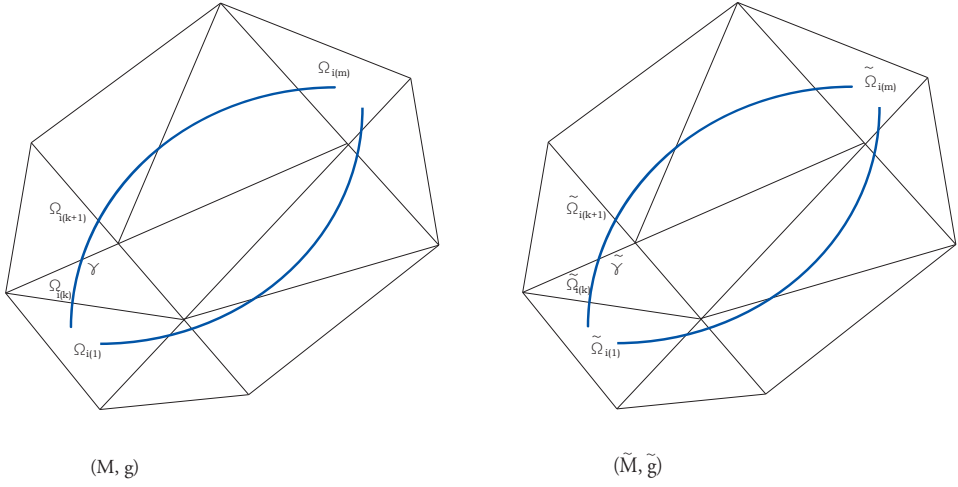}\\
  \caption{Chains of chambers, connecting points}\label{Fig3}
\end{figure}

We have proven that,
for any chain of $n-$simpleces $\Om_1=\Om_{i(1)},...,\Om_{i(m)}$ in
$\M$ which are pairwise adjoint, i.e. there is an $(n-1)$-interface
$\ga^{(k)}$ between $\Om_{i(k)}$ and $\Om_{i(k+1)},\,\,k=1,...,m-1,$
there is a chain of $n$-simpleces
$\wti\Om_1=\wti\Om_{i(1)},...,\wti\Om_{i(m)}$ in $\wti\M,$ such that
$\Om_{i(k)}$ and $\wti\Om_{i(k)}$ are diffeomorphic with
diffeomorphism $X_{i(k)}$ satisfying \mref{local} (strictly speaking
$X_{i(k)}$ may depend on a chain from $\Om_1$ to $\Om_{i(k)}$), see
Fig \ref{Fig3}.
Let us show that \begin{description}
    \item[(a)] For any $n-$simplex $\Om_j\subset\M$ there is an
    $n-$simplex $\wti\Om_j\subset\wti\M$ which is diffeomorphic to
    $\Om_j$ with diffeomorphism $X_j$ satisfying \mref{local} and,
    likewise, for any $\wti\Om_j\subset\wti\M$ there is $\Om_j
    \subset \M$ with described properties;
    \item[(b)] If $\Om_j$ and $\wti\Om_j$ are diffeomorphic with
    diffeomorphisms $X_j$ and $X_j^{'}$ satisfying \mref{local} then
    $X_j=X_j^{'};$
    \item[(c)] If $\Om_j$ is diffeomorphic to $\wti\Om_j$ and
    $\wti\Om_j^{'}$ by $X_j$ and $X_j^{'}$ satisfying \mref{local},
    then $\wti\Om_j = \wti\Om_j^{'}$.
\end{description} Let now $\Om_{j(1)},\,\wti\Om_{j(1)}$ be $\Om,\,
\wti\Om$ are the $n-$simpleces  with $\ga,\,\ti\ga$ used in Lemma
\ref{first} being their boundary $(n-1)-$subsimpleces. Then, by
Lemma \ref{first}, there is $X_1:\,\Om_{j(1)}\rightarrow\wti
\Om_{j(1)}$ satisfying \mref{local}. Denote by $\ga_1$ the
$(n-1)-$interface between $\Om_{j(1)}$ and $\Om_{j(2)}$. By
Corollary \ref{boundary?}, there is $\wti\Om_{j(2)}$ adjacent to
$\wti\Om_{j(1)}$ with the $(n-1)-$interface $\ti\ga_1=\wti\Om_{j(1)}
\cap\wti\Om_{j(2)}$ such that $\ti\ga_1=X_1(\ga_1)$. By Lemma
\ref{crossing}, there is an isometry $X_2:\,\Om_{j(2)}\rightarrow
\wti \Om_{j(2)}$ with $X_1|_{\ga_1}=X_2|_{\ga_1}$. Continuing this
process, we obtain an isometry $X_m:\Om_{j(m)}=\Om_j$, $\wti
\Om_{j(m)}:=\wti\Om_{j}$ which satisfies \mref{local}.

To prove (b), let $p\in\Om_j^{int}$, and $\ti p\in X_j(p)$, $\ti
p^{'}=X_j^{'}(p)\in\wti\Om_j^{int}$. As $X_j,\,X_j^{'}$ satisfy
\mref{local}, $\ti\varphi_k(\ti p)=\wti\varphi_k(\ti p^{'}),\,
k=1,2,\dots\,.$ By Proposition \ref{disting}, $\ti p=\ti p^{'}$. As
$p\in\Om_j^{int}$ is arbitrary, $X_j=X_j^{'}$ on $\Om_j$. Using the
fact that $\{\ti\varphi_k\}_{k=1}^\infty$ distinguish points in
$\wti\M^{int}$, see Proposition \ref{disting}, we prove property
(c).

Based on properties (a)-(c) we show the following result.
\begin{lemma}\label{almost}Let LBSD for the Riemannian polyhedron
$(\M,\,g)$ and $(\wti\M,\,\ti g)$ be equivalent. Then
\begin{enumerate}
    \item For any $\Om_i$ in $\M$ there is a unique diffeomorphism
    $X_i:\,\Om_i^{int}\rightarrow\wti\Om_i^{int}$, which satisfies
    \mref{local}. 
    \item $\ga$ is an $(n-1)-$interface between $\Om_-$ and $\Om_+$
    if and only if the corresponding $\ti\ga$ is an
    $(n-1)-$interface between $\wti\Om_-$ and $\wti\Om_+$. In this
    case $$X_-|_\ga=X_+|_\ga,$$ where $X_-$ is the closures of the
    described above diffeomorphisms $X_\mp$ on $\Om_\mp^{int}$,
    \item The diffeomorphisms $\cup X$ can be uniquely extended to
    an isometry $$X:\,\M^{reg}\rightarrow
    \wti\M^{reg}.$$ Here $X$ is an isometry of $\M^{reg}$ and
    $\wti\M^{reg}$ considered as metric space with the distance
    function given in Definition \ref{admissible}.
\end{enumerate} \end{lemma} \textbf{Proof} By (a)-(c) it remains to
prove the part of this lemma dealing with $(n-1)-$simpleces. Let
$\ga$ be an interface between $\Om_-$ and $\Om_+$ with $\wti\Om_-$
and $\wti\Om_+$ being the corresponding $n-$simpleces in $\wti\M.$
Crossing $\ga$ we move from $\wti\Om_-$ to $\wti\Om^{'}$ in $\wti\M$
so that, by the previous constructions, $\wti\Om^{int}$ is a
diffeomorphic to $\Om_+$ with diffeomorphism satisfying
\mref{local}. But also $\Om_+$ and $\wti\Om_+$ are diffeomorphic
with diffeomorphism satisfying \mref{local}. Thus $\wti\Om^{'}=
\wti\Om_+$ and the diffeomorphisms $X_{\mp}:\, \Om^{int}_\mp
\rightarrow \wti\Om_\mp^{int}$ are uniquely extendable to an
isometry $X_{-,+}: \,\Om^{int}_- \cup\Om_+^{int}\cup\ga\rightarrow
\wti\Om^{int}_- \cup\wti\Om_+^{int}\cup\ti\ga,$ where
$\ti\ga=X(\ga)$. As the distance functions in
$\M^{reg},\,\wti\M^{reg}$ employ only curves being in $\M\backslash
\M^{n-2},\,\wti\M\backslash\wti\M^{n-2}$, the above results
concerning the local isometries yield the desired global isometry.

Theorem \ref{main} immediately follows from Lemma \ref{almost}
taking into account the Definition \ref{admissible} of an admissible
Riemannian polyhedron. Note that this definition implies that the
topology of a Riemannian polyhedron considered as a metric space is
the same as the topology of the underlying polyhedron.

\section{Further generalizations and
estimates}\label{generalization}
\subsection{}Condition \ref{A} that the metric tensor does have a jump
singularity at every point of any $(n-1)-$dimensional interface,
i.e. Condition \ref{A} may be too restrictive. In this section we
relax it a bit. namely, we assume that it \begin{itemize}
    \item \emph{either} the metric tensor $g$ does have a jump singularity
    at any point $p\in\ga^{int},$ for an $(n-1)-$interface $\ga,$
    between $\Om_-$ and $\Om_+,$
    \item \emph{or} $g$ is smooth across $\ga$ at any $p\in\ga^{int}.$
\end{itemize} The latter condition means, that in the interface
coordinates related to $\ga,$ see section 1.2, the metric tensor
$g_{\al\be}(s,\si)$ is smooth at $\si=0.$ In this case we call $\ga$
an \emph{artificial} interface and would like to treat
$\Om_-\cup\Om_+$ together, introducing
$\Om^{int}_-\cup\ga^{int}\cup\Om_+^{int}.$ Further removing
artificial interfaces and taking into account that Riemannian
polyhedron consists of a finite number  of $n-$simpleces, it is
natural to introduce the following
object.\begin{definition}\label{chamber} A chamber $\Om$ of a
Riemannian polyhedron $(\M,\,g)$ is the maximal union of open
$n-$simpleces together with open artificial interfaces adjacent to
them.\end{definition} It is natural to treat each chamber $\Om$ as
an open Riemannian submanifold of $\M$ with a piece-wise smooth
boundary. However, due to \cite{Kervaire}, there may be topological
obstruction to that. Namely, it may happen that, for $p\in\M^{n-2},$
there is no open neighborhood of $p$ (in $\M$) which is homeomorphic
to $\mathbb{R}^n.$ In this connection we introduce
\begin{definition} \label{admissible1} An admissible Riemannian
polyhedron $(\M,g)$ is called weakly admissible if any of its open
chamber $\Om^{int}$ is an open $n-$dimensional Riemannian manifold
with piece-wise smooth boundary and, if $(n-1)-$subsimplex $\ga$ of
$\M$ is not artificial, they either $\ga\in\p\M$ or the metric $g$
has a jump singularity across $\ga.$  \end{definition} Then, similar
to Theorem \ref{main} we can prove
\begin{theorem}\label{main1} Let $(\M,\,g)$ and $(\wti\M,\,\ti g)$
be weakly admissible Riemannian polyhedra. Let local boundary
spectral data corresponding to $(\M,\,g)$ and $(\wti\M,\,\ti g)$ are
equivalent. Then there is an isometry $X:\,\M\rightarrow\wti\M$ such
that for any open chamber
$\Om^{int}_i\subset\M,\,\,X|_{\Om_i^{int}}$ is a diffeomorphism onto
an open chamber $\wti\Om_i^{int}\subset\wti\M,$ which satisfy
equation \mref{local}.
\end{theorem} Local boundary spectral data in this case consist of a
"smooth" open subset $\Ga_0\subset\p\M,$ such that there is an open
neighborhood $U\subset\M,\,$ of $\Ga_0,$ which is diffeomorphic to a
half-ball $\{x\in B(0,r):\,x_n\ge0\}$ with $\Ga_0$ being an open
subset of $B(0,r)\cap \{x_n=0\},$ and eigenpairs $\{\lambda_k,\,
\p_\nu\varphi_k|_{\Ga_0}\}_{k=1}^\infty.$

\subsection{} Similar to the smooth case, the uniqueness Theorem
\ref{main} remains valid when, instead of the boundary spectral
data, we have a local non-stationary Dirichlet-to-Neumann maps,
$$\Lambda_{\Ga_0}^T: \dot C^\infty(\Ga_0\times(0,T)) \rightarrow
\dot C^\infty (\Ga_0\times (0,T)),$$ where $\dot
C^\infty(\Ga_0\times(0,T))$ consists of smooth functions which are
qual to $0$ near $t=0,$ $\Lambda_{\Ga_0}^T$ is then defined as
$$\Lambda_{\Ga_0}^T: f\rightarrow \p_\nu u^f|_{\Ga_0 \times(0,T)},$$
where $u^f$ is the solution to the inverse boundary value problem
$$u^f_{tt}-\Delta u^f=0,\quad u^f|_{t<0}=0,\quad
u^f|_{\p\M\times(0,T)}=f.$$

Condition \ref{A} may be replaced by a condition that the metric
tensor $g$ is not smooth across any $(n-1)-$ interface $\ga$ at any
point $p\in\ga^{int}.$ This condition means that, for any
$p\in\ga^{int}$ and any interface conditions with $s=0$
corresponding to $p,$ the metric tensor $g_{\al\be}(0,\si)$ or some
its derivatives $\p_\si^k g_{\al\be} (0,\si)$ has a jump at $\si=0.$
Then Theorem \mref{main} remains valid as, for regions of $\ga$
where $g_{\al\be}(s,\si)$ is continuous across $\si=0$ together with
its first $(k-1)-$derivatives with respect to $\si,$ but $\p_\si^k
g_{\al\be}(0,\si),$ does have a jump across $\si=0,$ the reflected
Gaussian beam $U^r_\ve$ is of order $k.$ Namely, in representation
\mref{GB} for $U^{r}_\ve(t,y)$ the first non-zero $u_l(t,y)$ is
$u_k(t,y)$. This, however, does not alter considerations of section
\ref{GB}-\ref{globalRP} with the only difference that in section
\ref{51} the estimate for the reflected Gaussian beam in
\mref{furthoroL} is changed to $|\wti U _\ve^{r,L}(\ti q,
t)|>C\ve^k.$
 \subsection{Open problems}\begin{description}
 \item[(i)]This paper deals with uniqueness in the
 inverse boundary spectral problem for Riemannian polyhedron. It
 does not provide an algorithm to its reconstruction. In particular,
 trying to apply the technique of Section 4 \cite{KKL}, when
 approaching $\M^{n-2}$ the size of a step of the reconstruction
 procedure tends to zero. Therefore, it is impossible to reconstruct
 $(\M,\,g)$ by a finite number of steps.
 \item[(ii)] Even with generalization described in sections 5.1, 5.2,
 the class of Riemannian Polyhedron considered in this paper
 does not cover an important case when the metric $g$ is smooth
 across some part of the interface $\ga$ but is not smooth across
 other part of $\ga,$ i.e. we have "holes" in interfaces. We intend
 to study such cases in the forthcoming paper.
 \item[(iii)]Another important open question is the one of stability
 which could relate observation error with , in addition to
 curvature, injectivity radii, etc., \cite{AKKLT}, size of
 $n-$simpleces, value of jumps, etc.
 \end{description}\section {Acknowledgements}

The authors would like to express their gratitude to Prof. M. Lassas
for numerous stimulating discussions regarding the various aspects
of the problem, Prof. Yu. Burago and Dr. N. Kossovski for
consultations on the geometric issues and Prof. A. Kachalov for the
useful discussions on the asymptotical aspects of non-stationary
Gaussian beams.

The research of the first author was financially supported EPSRC
grant Ep/D065771/1.


\begin{thebibliography}{99}

\addcontentsline{toc}{section}{\texttt{References}}
%

\bibitem{AMR} Alessandrini, G., Morassi, A., Rosset, E. Detecting cavities
by electrostatic boundary measurements, {\it Inv. Probl.,} {\bf 18} (2002),  1333--1353.
%
%
%


\bibitem{PaiAst} K. Astala, M. Lassas, L. P\"aiv\"arinta, Calderon's
inverse problem for anisotropic conductivity in the plane,
Communications in Partial Differential Equations \textbf{30}(2005),
no. 1-3, 207-224.

\bibitem{ALP} Astala, K., P\"aiv\"arinta, L., Lassas, M. Calder—n's
inverse problem for anisotropic conductivity in the plane, {\it
Comm. Part. Diff. Eq.}, {\bf 30} (2005),  207--224.

\bibitem{AKKLT} M. Anderson, A.Katsuda, Ya. Kurylev, M. Lassas, M. Taylor,
Boundary Regularity for the Riccati equation, Geometry Convergence
and Gelfand's Inverse Boundary Problem, Inventiones Mathematicae
\textbf{158}(2004), 261-321.

\bibitem{BabichUlin} V. Babich, V. Ulin, The complex space-time
ray method and "quasiphotons", (Russian) Zap. Nauchn. Sem. LOMI,
\textbf{117}(1981), 5-12.

%
\bibitem{WBallman} W. Ballman, A volume estimate for
piecewise smooth metrics on simplicial complexes, Rendiconti Sem.
Mat. Fis. Milano \textbf{66}(1996), 323-331.

\bibitem{Bel}
M. Belishev, An approach to multidimensional inverse problems for
the wave equation, (Russian) Dokl. Akad. Nauk SSSR
\textbf{297}(1987), no.3, 524-527; translated in Soviet Math. Dokl.
\textbf{36}(1988),no.3 481-484.

\bibitem{Belishev} M. Belishev, Wave basis in multidimensional
inverse problems, (Russian) Mat.Sb. \textbf{180}(1989), 584-602.

%

\bibitem{BelKur}
M. Belishev, Ya. Kurylev, To the construction of a Riemannian
manifold via its spectral data (BC-method), Communications in
Partial Differential Equations \textbf{17}(1992), no.5-6, 591-594.
%
%
%
%
%
%
%
%
%
%
%

\bibitem{FugEells} J. Eells, B. Fuglede, \emph{Harmonic maps between
Riemannian polyhedra.} With a preface of Gromov M., Cambridge Tracts
in Mathematics, \textbf{142}, Cambridge University Press, Cambridge,
2001.
%
%
%
%
%
\bibitem{GrUhl} Greenleaf, A.,  Uhlmann, G. Recovering singularities
of a potential from singularities of scattering data, {\it Comm.
Math. Phys.}, {\bf 157} (1993),  549--572.
%
%
%

\bibitem{Federer}
H. Federer, \emph{Geometric Measure Theory,} Springer-Verlag,
Berlin, Heidelberg, New York, 1996.

\bibitem{Fuglede} B. Fuglede, Finite energy maps from
Riemannian polyhedra to metric spaces, Ann.Acad.Sci.Fenn.Math.,
\textbf{28}(2003), no.2, 433-458.


\bibitem{Ike} Ikehata, M. Reconstruction of inclusion from boundary measurements,
{\it  J. Inverse Ill-Posed Probl.}, {\bf 10} (2002),  37--65.

\bibitem{Ike1}  Ikehata, M.,  Siltanen, S. Electrical impedance tomography and
Mittag-Leffler's function, {\it  Inv. Prob.}, {\bf 20} (2004), 1325--1348.

\bibitem{Isakov} Isakov V., On uniqueness of recovery of a discontinuous
conductivity coefficient, Comm.Pure Appl.Math., \textbf{41}(1988),
no.7, 865-877.

\bibitem{Isakov1} Isakov V.  Inverse Problems for Partial Differential Equations. Springer, New York, 2006, 344 pp.

\bibitem{Iso}
Isozaki, H. Inverse scattering theory for Dirac operators.
{\it Ann. Inst. H. Poincare Phys. Theor.}, {\bf 66} (1997),  237-270.


\bibitem{Katchalov}
A. Kachalov, Gaussian beams, Hamilton-Jacobi equations and Finsler
geometry, Zapiski Nauchn. Semin.POMI, \textbf{297}(2003), 66-92.


\bibitem{KK2}
Kachalov, A., Kurylev, Y.
Multidimensional inverse problem with incomplete boundary
spectral data.
{\it Comm. PDE}, {\bf 23} (1998), 55-95.


\bibitem{KKL}
A. Kachalov, Ya. Kurylev, M. Lassas,\emph{ Inverse Boundary Spectral
Problems}, Chapman Hall / CRC {\bf 123}, 2001.

\bibitem{KKL2} A. Katchalov, Y. Kurylev, M. Lassas, \emph{ Energy measurements and
equivalence of boundary data for inverse problems on non-compact
manifolds}. IMA volumes in Mathematics and Applications (Springer
Verlag) ``Geometric Methods in Inverse Problems and PDE Control''
Ed. C. Croke, I. Lasiecka, G. Uhlmann, M. Vogelius, (2004), 183-214.
%

\bibitem{KSU}
Kenig C., Sjoestrand  J., Uhlmann  G.
The Calder\'on problem with partial data, preprint
arXiv math.AP/0405486.

\bibitem{Kervaire} M. Kervaire, A manifold which does not
admit any differentiable structure, New York (USA), Commentarii
mathematici Helvetici, \textbf{34}(1960), 257-270.

\bibitem{KirpichGb} A. Kirpichnikova, Propagation of a Gaussian
beam near an interface  in an anisotropic medium, Zapiski Nauchn.
Sem. POMI, \textbf{324(34)}(2005), 77-109.

\bibitem{KirPai}
A. Kirsch, L. P\"aiv\"arinta, On recovering obstacles inside
inhomogeneities, Math.Meth.Appl.Sci., \textbf{21}(1998), 619-651.
%
%
%

\bibitem{KLS} Kurylev y., Lassas M., Somersalo E.
Maxwell's equations with a polarization independent wave velocity:
Direct and inverse problems,
{\it J.  Mathem. Pures  Appl.}
{\bf 86} (2006), 237-270.

\bibitem{KL-Dirac} Kurylev Y., Lassas M.
Inverse problem for a Dirac-type equation on a vector bundle,
arXiv math.AP/0501049.

%

\bibitem{LTU}
Lassas, M. Taylor, M., Uhlmann, G.
The Dirichlet-to-Neumann map for complete Riemannian manifolds with
boundary,
{\it Comm. Anal. Geom.} {\bf 11} (2003), 207-22.

\bibitem{LeU}
Lee, J., Uhlmann, G.
Determining anisotropic real-analytic conductivities by boundary
measurements.
{\it Comm. Pure Appl. Math.} {\bf 42} (1989), no. 8, 1097--1112.

\bibitem{Na1}
Nachman, A.
    Reconstructions from boundary measurements.
{\it  Ann. of Math.} (2) {\bf
128} (1988),
no. 3, 531--576.

\bibitem{Na2}
Nachman, A. Global uniqueness for a two-dimensional
    inverse boundary value problem. {\it Ann.
of Math.} (2) {\bf 143} (1996), no. 1, 71--96.

\bibitem{NakUhl}Nakamura, G., Uhlmann, G.  Global uniqueness for an
 inverse boundary value problem arising in elasticity,  {\it Invent. Math.} {\bf 118} (1994),
 457--474, {\bf 152} (2003), 205--207.

\bibitem{Nak-Dirac} Nakamura, G.,  Tsuchida, T. Uniqueness for an
inverse boundary value problem for Dirac operators. Comm. Part.
Diff. Eq. {\bf 25} (2000), 1327--1369.

\bibitem{Nov} Novikov R. A multidimensional inverse spectral problem for the
equation: $\Delta \psi+(v(x)-Eu(x))\psi=0$, Funkt. Anal. i ego Priloz. (in Russian),
{\bf 22} (1988), 11-22.

\bibitem{OPS} Ola, P.,  P\"aiv\"arinta, L.,  Somersalo, E. An inverse boundary
value problem in electrodynamics,  {\it Duke Math. J.}, {\bf 70}
(1993),  617--653.

\bibitem{OS} Ola, P., Somersalo, E. Electromagnetic inverse problems and
generalized Sommerfeld potentials, {\it SIAM J. Appl. Math.},
{\bf 56} (1996), 1129--1145.


\bibitem{PU}
 Pestov, L., Uhlmann, G. Two dimensional compact simple Riemannian manifolds are boundary distance rigid.  {\it Ann. of Math.}, {\bf 161}  (2005),
1093-1110.


\bibitem{PopovMM} M. Popov, \emph{Ray Theory and Gaussian Beam Method
for Geophysicists}, Edufba, Salvator-Bahia, 2002.

\bibitem{Ralston} J. Ralston, \emph{Gaussian beams and propagation
of singularities,} Studies in PDE, MAA Studies in Mathematics,
\textbf{23}, Walter Littman ed, 1983.
%
%

\bibitem{Syl}
Sylvester, J. An anisotropic inverse boundary value problem.
{\it Comm. Pure Appl. Math.} 43 (1990), no. 2, 201--232

%
%
%
%
%

\bibitem{SU} Sylvester J., Uhlmann G.
A global uniqueness theorem for an inverse
boundary value
problem. {\it Ann. of Math.}, {\bf 125} (1987),  153-169.
%
%
%

\bibitem{Taylor} M. Taylor, \emph{Tools for PDE}, AMS, Providence,
R.I., 2000.

%
%
\end{thebibliography}
\end{document}